\documentclass{amsart}
\usepackage[latin1]{inputenc}
\usepackage{amsmath,amsthm,amsfonts,latexsym,amscd,amssymb,enumerate}

\DeclareMathOperator{\Tr}{Tr}
\DeclareMathOperator{\TR}{TR}
\DeclareMathOperator{\tr}{tr}
\DeclareMathOperator{\res}{res}
\DeclareMathOperator{\op}{op}
\DeclareMathOperator{\ord}{ord}

\DeclareMathOperator{\ind}{Index}

\newcommand{\C}{{\mathbb C}}

\newcommand{\R}{{\mathbb R}}

\newcommand{\Z}{{\mathbb Z}}

\newcommand{\cD}{{\mathcal D}}

\newcommand{\go}{\omega}

\newcommand{\gvp}{\varphi}

\newcommand{\gs}{\sigma}

\newcommand{\forget}[1]{}

\newcommand{\psido}{$\psi$do}
\renewcommand{\Re}{{\rm Re}\,}

\newcommand{\bli}[1]{\begin{list}{{\rm(#1{num})}\hfill}{\usecounter{num}\labelwidth1cm
\leftmargin1cm\labelsep0cm\rightmargin1pt\parsep0.5ex plus0.2ex minus0.1ex
\itemsep0ex plus0.2ex\itemindent0cm}}
\newcommand{\eli}{\end{list}}

\title[Uniqueness of the Kontsevich-Vishik Trace]%
{Uniqueness of the Kontsevich-Vishik Trace}
\author{L.\ Maniccia}
\address{Universit\`a di Bologna, Dipartimento di Matematica, Piazza di Porta S. Donato 5, 
         40127 Bologna, Italy}
\email{maniccia@dm.unibo.it}

\author{E.\ Schrohe}
\address{Leibniz Universit\"at Hannover, Institut f\"ur Analysis, Welfengarten 1, 
30167 Hannover, Germany}
\email{schrohe@math.uni-hannover.de}

\author{J.\ Seiler}
\address{Leibniz Universit\"at Hannover, Institut f\"ur Angewandte Mathematik, \mbox{Welfengarten 1}, 
30167 Hannover, Germany}
\email{seiler@ifam.uni-hannover.de}

\begin{document}
\begin{abstract} Let $M$ be a closed manifold. 
We show that the Kontsevich-Vishik trace, which is defined on the set 
of all classical pseudodifferential operators
on $M$, whose (complex) order is not an integer 
greater than or equal to $- \dim M$, 
is the unique functional which (i) is linear on its domain, 
(ii) has the trace property and (iii) coincides with the $L^2$-operator 
trace on trace class operators. 

Also the extension to 
even-even pseudodifferential operators of arbitrary integer order 
on odd-dimensional manifolds and to even-odd  pseudodifferential operators 
of arbitrary integer order on even-dimensional manifolds is unique.  
\end{abstract}

\maketitle

\section{Introduction}
We denote by $M$ a compact $n$-dimensional manifold without boundary.
A classical pseudodifferential operator  ($\psi$do) 
$A$ on $M$ is said to have order $\mu\in\C$, if it belongs to the Hörmander 
class $S_{1,0}^{\,\Re \mu}(M)$ and the local symbols $a=a(x,\xi)$ of $A$ have 
asymptotic expansions 
\begin{eqnarray}
\label{asexp}
a\sim\sum_{j=0}^\infty a_{\mu-j},
\end{eqnarray} 
where the $a_{\mu-j}$ are positively 
homogeneous of degree $\mu-j$ for large $\xi$.
In general, all pseudodifferential operators will be assumed to act on sections 
of vector bundles over $M$.
We shall write $\ord A=\mu$ to express that the order of $A$ is $\mu$. 

In two remarkable papers, Kontsevich and Vishik in 1994 and 1995 
analyzed the properties of determinants of elliptic \psido's, 
\cite{KoVi1}, \cite{KoVi2}.
One important tool was the construction of a trace-like mapping $\TR$ 
defined on the set of all classical \psido's whose order is not an 
element of $\Z_{\ge -n}$, the set of integers greater than or equal to $ -n$. 
 
We shall denote this domain by $\cD$. As the sum of two operator of orders 
$\mu$ and $\mu'$ in $\cD$ is an element of $\cD$ only if $\mu- \mu'$ is an integer,
$\cD$ is not a vector space. 
Thus it does not make sense to speak about linear functionals on $\cD$.
The map  $\TR:\cD\to \C$, however, is as linear as it 
can be expected to be: 
\begin{eqnarray}
\label{tr1}
\mbox{\ \ }\TR(cA+dB)=c\TR(A) + d\Tr(B)\quad\text{for } c,d\in \C, A,B,cA+dB\in \cD.
\end{eqnarray}
Moreover, $\TR$ behaves like a trace:
\begin{eqnarray}
\label{tr2}
\TR(AB)=\TR(BA),\quad\text{whenever }AB, BA\in \cD.
\end{eqnarray}
Finally, the Kontsevich-Vishik trace (sometimes also: canonical trace) 
$\TR$ coincides with the $L^2$-operator trace on
trace class \psido's:
\begin{eqnarray}
\label{tr3}
\TR(A)=\tr(A)\quad \text{if }\Re\ord(A)<-n.
\end{eqnarray} 
It is clear that the Kontsevich-Vishik trace cannot be extended to a trace
on the algebra of all \psido's on $M$: 
The only trace there (up to multiples) 
is the Wodzicki residue \cite{W}, which is known to vanish on trace class operators. 
There also is a simple direct way to see this: 
We know -- e.g. from the Atiyah-Singer index theorem -- that 
there exists an elliptic pseudodifferential operator $P$  on $M$ with nonzero index.
Using order reducing operators, we may assume the order of $P$ to be 
zero. 
Let  $Q$ be a parametrix to $P$ modulo smoothing operators. Then 
$$\ind P = \tr(1-PQ)-\tr(1-QP).$$
If we could extend $\TR$ to a trace on all pseudodifferential operators, the 
right hand side 
could be rewritten as the trace of the commutator $\TR [P,Q]$ 
and therefore would have to be zero -- a contradiction. 

It has been observed, however, by Kontsevich-Vishik and Grubb \cite{Grubb} 
that $\TR$ extends to a slightly larger domain.
Recall that the symbol $a$ of an integer order operator $A$ is said to be 
even-even, if the homogeneous components satisfy 
\begin{eqnarray}
\label{evev}
a_{\mu-j}(x,-\xi)=(-1)^{\mu-j} a_{\mu-j}(x,\xi).
\end{eqnarray}
It is called even-odd, if 
\begin{eqnarray}
\label{eveod}
a_{\mu-j}(x,-\xi)=(-1)^{\mu-j+1} a_{\mu-j}(x,\xi).
\end{eqnarray}
The Kontsevich-Vishik trace $\TR(A)$ for a \psido \ $A$ of order $\mu$ 
then can also be defined if $\mu\in\Z_{\ge-n}$,
provided that
\begin{itemize}
\item[(EE)] 
 $n$ is odd, and the symbol of $A$ is even-even, or
\item [(EO)]
 $n$ is even, and the symbol of $A$ is even-odd.  
\end{itemize}
For the sake of brevity we shall denote this larger domain (depending on $n$) by $\cD^+$.

In both cases, the component $a_{-n}$ in the asymptotic expansion of the 
symbol of $A$ is odd in $\xi$ for large $|\xi|$, say for $|\xi|\ge 1$:
$$ a_{-n}(x,-\xi)=-a_{-n}(x,\xi).$$ 
Hence the density for the Wodzicki residue of the operator $A$ vanishes
pointwise, i.e. 
$$\res_x(A)=\int_{S^*_xM} a_{-n}(x,\xi)\,d\gs(\xi)=0\quad\text{ for each }x\in M. $$
Here, $d\gs$ is the surface measure on the unit sphere $S_x^*M$ over $x$ in the cotangent 
bundle. The Wodzicki residue of $A$ is given by integration of $\res_x A$ over $M$
and therefore also vanishes.

The trace property \eqref{tr2} extends to the case where $A$ and $B$ have integer order 
and $AB$ and $BA$ belong to $\cD^+$.

The Kontsevich-Vishik trace has received considerable attention and found interesting 
applications, see e.g. \cite{Grubb2,Lesc,MSS,Okik,PS}.
Moreover, it has been extended to boundary value problems in 
Boutet de Monvel's calculus \cite{GS2}.

It seems, however, that it never has been noticed that 
the above properties make the Kontsevich-Vishik trace 
unique. 
This is what we show in this short note. The proof, which will be given in the next section, 
relies on ideas in \cite{FGLS}.

\medskip\noindent{\bf Theorem.}
{\bli\alph\item Let $\tau:\cD\to \C$ be a map with properties 
\eqref{tr1}, \eqref{tr2}, and \eqref{tr3}. Then $\tau=\TR$.
\item Also the extension of $\tau$ to $\cD^+$ is unique. In fact, $\tau$ already 
is unique on the space of all integer order \psido's which satisfy {\rm (EE)} or 
{\rm (EO)} when $\mu\ge -n$.
\eli}

\section{Proof}

In order to establish (a),  
choose a \psido\ $A$ of order $\mu\in\C\setminus \Z_{\ge -n}$ on $M$. 

We find a covering of $M$ by open neighborhoods and a finite subordinate partition 
of unity $\{\gvp_j\}$  such that for every pair $(j,k)$, both $\gvp_j$ and $\gvp_k$ have
support in one coordinate neighborhood. 
We write 
$$A=\sum_{j,k}\gvp_j A\gvp_k.$$
Each operator $\gvp_j A\gvp_k$ may be considered a \psido\ on $\R^n$.
As the map $\tau$ has the linearity property \eqref{tr1}, 
we may confine ourselves to the case where $A=\op a$ with a  
symbol $a$ on $\R^n$ having an expansion \eqref{asexp}. Moreover, we can assume
that $A=\gvp A\psi$ whenever $\gvp,\psi\in C^\infty_c(\R^n)$ are equal to one on
a sufficiently large set.

To simplify further, we write 
\begin{eqnarray}
\label{op}
A=\op a_\mu + \op a_{\mu-1} +  \ldots \op a_{\mu-K} +\op r,
\end{eqnarray} 
where $a_{\mu-j}$ is a symbol on $\R^n$, homogeneous in $\xi$  of degree $\mu-j$
for $|\xi|\ge 1$, and $K$ is so large that 
$r\in S^{-n-1}_{1,0}$.
For $\gvp,\psi\in C^\infty_c(\R^n)$ as above we then have 
$$\tau(\op a)=\tau(\gvp\op (a)\psi) = 
\sum_{j=0}^K\tau(\gvp\op (a_{\mu-j})\psi)+\tau(\gvp\op (r)\psi).$$
Since $\tau(\gvp\op( r)\psi)=\tr(\gvp\op( r)\psi)$ by \eqref{tr3}, 
we will know $\tau(\op a)$ as soon as 
we know $\tau(\gvp\op(a_{\mu-j})\psi)$ for $j=0,\ldots,K$.

We may assume that $\mu$ is not an integer, since the operator trace determines 
$\tau$ on all operators of order  $\mu<-n$. 
Now we let 
\begin{eqnarray}
\label{b}
b_{\mu-j}(x,\xi)
=\frac1{n+\mu-j}\ \sum_{k=1}^n \partial_{\xi_k}(\xi_k a_{\mu-j}(x,\xi)).
\end{eqnarray}
Euler's relation for homogenous functions implies that, for $|\xi|\ge 1$, 
\begin{eqnarray*}
b_{\mu-j}=\frac1{n+\mu-j}\left(n a_{\mu-j}+(\mu-j)a_{\mu-j}\right)=a_{\mu-j}.
\end{eqnarray*}
Hence we can write 
\begin{eqnarray}
\label{split}
\tau(\gvp\op( a_{\mu-j})\psi)= 
\tau(\gvp\op (a_{\mu-j}- b_{\mu-j})\psi)+ \tau(\gvp\op( b_{\mu-j})\psi). 
\end{eqnarray}   
Since $a_{\mu-j}- b_{\mu-j}$ is regularizing, the first term on the right hand side
is determined by Property \eqref{tr3}.
Now we choose additionally  
$\chi\in C^\infty_c(\R^n)$ with $\chi\gvp=\gvp$ and $\chi\psi=\psi$. 
The fact that 
$\op(\partial_{\xi_k}p)=-i\ [x_k,\op p\,]$ for an arbitrary symbol $p$ implies that 
$$\gvp \op(b_{\mu-j})\psi=-i \sum_{k=1}^n [\chi x_k,\gvp\op(\xi_ka_{\mu-j})\psi].$$
   
Assuming that $\tau$ has Property \eqref{tr2}, it vanishes on the 
last term in \eqref{split} which 
is a sum of commutators.
Hence the proof of (a) is complete.

Next let us show (b). 
With the same considerations 
as before we may assume that $A=\op a$ 
is a pseudodifferential operator on $\R^n$ with a representation 
as in \eqref{op}, where now $\mu$ is an integer $\ge -n$ and the $a_{\mu-j}$
have property (EE) or (EO).
We only have to show that $\tau(\gvp\op(a_{\mu-j})\psi)$ is uniquely determined,
$j=0,\ldots,\mu+n$. For $\mu-j\not=-n$ the argument is as before, using the symbols 
in \eqref{b} and noting that 
$a_{\mu-j}(x,\xi)\xi_k$ is even-even or even-odd whenever 
this is the case for $a_{\mu-j}$.
 
So let us consider $a_{-n}$.  
The assumption that $n$ is odd and 
$a_{-n}$ even-even or $n$ is even and 
$a_{-n}$ even-odd implies that $a_{-n}$ is odd in $\xi$:
$$a_{-n}(x,-\xi)=-a_{-n}(x,\xi)\quad\text{for }|\xi|\ge 1.$$
Hence, for each fixed $x$, the integral over the unit sphere $S=\{|\xi|=1\}$ vanishes:
\begin{eqnarray}
\label{int=0}
\int_Sa_{-n}(x,\xi)\, d\gs(\xi)=0.
\end{eqnarray}

The Laplace operator $\Delta=\sum_{k=1}^n{\partial^2}/{\partial \xi_k^2}$ 
in polar coordinates  takes the form  
$$\Delta = \frac1{r^{n-1}}\frac\partial{\partial r}\left(r^{n-1}\frac\partial{\partial r}\right)
+\frac1{r^2}\Delta_S,$$
where $r=|\xi|$ is the radial variable and $\Delta_S$ is the Laplace-Beltrami 
operator on $S$.

Equation \eqref{int=0} implies that, for each $x$,  the function 
$a_{-n}(x,\cdot)$ is orthogonal to the constants which form the kernel of
the symmetric operator $\Delta_S$.
Hence there is a unique function $q(x,\cdot)\in C^\infty(S)$, orthogonal
to the constants, such that 
$\Delta_Sq(x,\cdot)=a_{-n}(x,\cdot)|_S$.
As $\Delta_S$ commutes with the antipodal map
$\eta\mapsto -\eta$, we have 
$\Delta_S(q(x,-\cdot))=a_{-n}(x,-\cdot)|_S=-a_{-n}(x,\cdot)|_S$.
Hence $q(x,\cdot)+q(x,-\cdot)$ belongs to the kernel of $\Delta_S$,
thus is constant. On the  other hand, both  $q(x,\cdot)$ and $q(x,-\cdot)$
are orthogonal to the constants. Therefore $q(x,\cdot)+q(x,-\cdot)$ is zero, 
i.e., $q(x,\cdot)$ is an odd function on $S$.   

Now we choose a smooth function $\go$ on $\R$ which vanishes for small $r$ 
and is equal to $1$ for $r\ge 1/2$. We  let
$$b_{-n}=\go(r)r^{2-n}q=\go(|\xi|)\,|\xi|^{2-n}\,q(x,\xi/|\xi|).$$
This is a smooth function on $\R^n$ which is homogeneous of degree $2-n$ in 
$\xi$ for $|\xi|\ge 1$.
As $a_{-n}(x,\xi)$ vanishes for $x$ outside a compact set, so does $b_{-n}(x,\xi)$.
In particular, $b_{-n}$ is an element of 
$S^{2-n}_{1,0}(\R^n\times\R^n)$.
Moreover, we have for $|\xi|\ge 1$ 
\begin{eqnarray*}
\Delta b_{-n}=\Delta(r^{2-n}q(x,\cdot))=r^{-n}a_{-n}(x,\cdot)|_S= a_{-n}.
\end{eqnarray*}
We write $a_{-n}=(a_{-n}-\Delta b_{-n})+\Delta b_{-n}$. 
The symbol $a_{-n}-\Delta b_{-n}$ is regularizing and  
thus $\tau(\gvp\op(a_{-n}-\Delta b_{-n})\psi)$ is determined by \eqref{tr3}.
The operator associated with $\op(\gvp(\Delta b_{-n})\psi)$ 
on the other hand is a sum of commutators:
\begin{eqnarray}
\label{a}
\gvp\op(\Delta b_{-n})\psi=-i\sum_{k=1}^n\ [\chi x_k,\gvp\op(\partial_{\xi_k}b_{-n})\psi],  
\end{eqnarray}
where $\chi$ is chosen as in the proof of (a).
Hence $\tau$ vanishes on $\gvp\op(\Delta b_{-n})\psi$. 
This concludes the argument.
 
\begin{small}
\bibliographystyle{amsalpha}

\end{small}
\end{document}